\documentclass[a4paper,10pt]{siamltex}
\usepackage{a4wide}

\usepackage{here}
\usepackage{amsmath}
\usepackage{amssymb}
\usepackage{amsfonts}
\usepackage[margin=1.5cm]{caption}
\usepackage[utf8]{inputenc}
\usepackage{lmodern}
\usepackage[english]{babel}
\usepackage{graphicx,multirow}
\usepackage{subfigure}

\newcommand*\samethanks[1][\value{footnote}]{\footnotemark[#1]}

\def\E{{\mathbb E}}
\def\R{{\mathbb R}}
\emergencystretch=\hsize
\begin{document}

\title{Adaptive numerical integration and control variates for pricing
Basket Options} 
\author{Christophe De Luigi\thanks{Aix Marseille Universit\'e, CNRS, ENSAM, LSIS, UMR 7296, 13397 Marseille, France; Universit\'e de Toulon, CNRS, LSIS, UMR 7296, 83957 La Garde, France; \{deluigi,maire\}@univ-tln.fr} \and J\'er\^ome
Lelong\thanks{Grenoble INP, Laboratoire Jean Kuntzmann, CNRS UMR 5224, 51 rue des
Mathématiques, 38041 Grenoble cedex 9, France; jerome.lelong@imag.fr} \and
Sylvain Maire\samethanks[1]} \maketitle

\begin{abstract}
  We develop a numerical method for pricing multidimensional vanilla options
  in the Black-Scholes framework. In low dimensions, we improve an adaptive
  integration algorithm proposed by two of the authors by introducing a new
  splitting strategy based on a geometrical criterion. In higher dimensions,
  this new algorithm is used as a control variate after a dimension reduction
  based on principal component analysis. Numerical tests are performed on the
  pricing of basket, put on minimum and digital options in dimensions up to
  ten.  
\end{abstract}

\section{Introduction}

Pricing financial derivatives generally boils down to numerically compute an
expectation, at least for European type contracts. Basically, there are two
ways of doing so: either solving a partial differential equation or resorting
to Monte Carlo techniques. Monte Carlo integration is known to provide better
results when the dimension of the problem increases but for reasonably small
dimensional problems, its efficiency is not that clear.  However, alternative
numerical integration techniques may help in such situations.  In the
Black-Scholes framework, the pricing of vanilla options reduces to a
numerical integration problem over $\mathbb{R}^{d}$, where $d$ is the number
of underlying assets. This problem has some specificities coming from the
properties of the function to integrate.  First, the integrand decreases
quickly away from the origin and hence we can consider that the integration
domain is a hypercube $[-A,A]^{d}$ with $A$ relatively small.  Second, the
integrand is clearly not a smooth function in the whole domain as it is only
continuous at the interface between the area where the function vanishes and
the area where it is positive. Moreover, this interface is neither known and
nor located at the boundary of the hypercube so there is no hope of using
techniques like the periodisation method \cite{key-HELLUY} to increase the
smoothness of the integrand. Finally, the integrand is a muldimensional
function in a dimension $d$, which can be large, so one also has to deal with
the curse of dimensionality. 

Hence, it seems quite natural to use Monte Carlo or quasi-Monte Carlo
\cite{key-NIEDERREITER} methods to face such a difficult problem. In fact a
crude use of these methods is not necessarily sufficient to reach a good
accuracy. Besides the usual variance reduction methods, one needs to develop
adaptive methods to make them really competitive. In many situations, the
function to integrate or to approximate may have completely different
behaviors in terms of variations or even in terms of regularity in different
parts of the domain $D.$ In those cases, it might be more efficient to
adaptively split $D$ in subregions according to error indicators based on
quadrature points. The most famous adaptive Monte Carlo integration routines
are MISER \cite{key-PRESS} and VEGAS \cite{key-LEPAGE}. They rely on
stratified and importance sampling and error indicators based on the
empirical variance, respectively. Quasi-Monte Carlo versions of these two
algorithms have been introduced in \cite{key-SCHURER1}. More recently,
adaptive approaches have been developed for stratified sampling
\cite{key-ETORE,key-JOURDAIN} and for importance sampling \cite{key-LAPEYRE}.
The idea of adaptive algorithms is to make use of past simulations to help
better leading future ones. In the case of stratified sampling, both the
number of samples in each strata and the boundaries of the strata can be
learnt inline. The idea remains quite similar for adaptive importance
sampling, which consists in adaptively learning the optimal change of measure
from the already drawn samples.  Following the methodology of
\cite{key-SCHURER2}, two of the authors have developed an adaptive
integration algorithm \cite{key-CDL} based on quasi-Monte Carlo quadratures,
which proved to be efficient for very smooth but also for less smooth
functions \cite{key-MAIRE4}.  This algorithm has already been tested on the
pricing of basket options achieving excellent results in dimension two but
loosing most of its efficiency for larger dimensions. The loss of accuracy is
mainly due to the difficulty to capture the interface between the two regions
of interest.

Our goal is to improve this algorithm in order to use it in higher
dimensions, for the pricing of more general options like digital or
put on minimum options and also for the computations of sensitivities with
respect to the parameters of the Black-Scholes model. The improvements
of the algorithm will consist in a new splitting criterion and on
a dimension reduction using a principal components analysis combined with
control variates. 

The rest of the paper is organized as follows. In section 2, we describe the
different kinds of options considered in the Black-Scholes framework and we
discuss the type of sensitivities we are interested in. In section 3, we
recall the adaptive algorithm developed in \cite{key-CDL} and introduce a new
criterion based on very simple geometric considerations like the ones used in
mesh refinement for finite element methods. In section 4, we compare the new
criterion to the old one on the pricing of different types of options in low
dimensions. We focus on examples in dimension two in order to emphasize the
quality of the mesh refinement. Section 5 is devoted to sensitivity analysis.
We compute the Delta of the different derivatives using standard
deterministic techniques like polynomial interpolation allowed by the high
accuracy of our adaptive pricing method. The last section deals with higher
dimensional models. The adaptive method is coupled with dimension reduction
via a principal component analysis to derive a variance reduction method based
on control variates.

\section{Presentation of the model}

\subsection{The model framework}
\label{sec_model}

We consider a Black-Scholes model in dimension $d$ in which each asset is
supposed to follow the standard one dimensional dynamics given under the risk
neutral measure by 
\[
dS_{t}^{i}=S_{t}^{i}(rdt+\sigma_{i}dW_{t}^{i})
\]
with $S_{0}^{i}=s^{i}$ and where $W_{t}=(W_{1}^{t},...,W_{d}^{t})$ denotes a
vector of correlated standard Brownian motions.  The volatility $\sigma$ is a
vector in $\mathbb{R}^{d}$, the instantaneous interest rate is $r$ and
$(s^{1},...,s^{d})$ is the vector of spot values. The covariance structure of
these correlated Brownian motions is supposed to be defined by $\left\langle
W,W\right\rangle _{t}=\Gamma t$ where $\Gamma$ is a positive definite matrix
with all its diagonal terms equal to one. In the numerical examples
considered in the next sections, we assume that 
\[
\Gamma_{i,j}=\delta_{i,j}+\rho(1-\delta_{i,j})
\]
where the parameter $\rho\in]-\frac{1}{d-1},1[$ to ensure that the matrix
$\Gamma$ remains positive definite. We introduce the Cholesky decomposition
$C$ of $\Gamma$ (such that $CC^{T}=\Gamma)$ and denote  by
$C_{i}$ its $i^{th}$ row for $1\leq i\leq  d$. The Black-Scholes model can be rewritten 
\[
S_{t}^{i}=S_{0}^{i}\exp\left\{(r-\frac{\sigma_{i}^{2}}{2})t+\sigma_{i}C_{i}B_{t}\right\}
\]
where now $B$ is a standard $d-$dimensional Brownian motion.

In this model, we want to price options with payoffs written
as functions of the asset price at a maturity time $T$. Hence, the price
is given by the discounted expectation $\exp(-rT)\E(\psi(S_{T}))$
where the function $\psi$ characterizes the option type. In the following,
we consider three different multidimensional options for
\begin{itemize}
  \item Basket options  with payoffs
    \[
    \psi(S_{T})=\left(\sum_{i=1}^{d}\lambda_{i}S_{t}^{i}-K\right)_{+},
    \]
  \item  Digital basket options with payoffs
    \[
    \psi(S_{T})=\left(\sum_{i=1}^{d}\lambda_{i}S_{t}^{i}-K\right)_{+}
    1_{\{\forall\,1\leq i\leq d,\, S_{t}^{i}\leq U_{i}\}},
    \]
    
  \item Put on minimum options with payoffs
    \[
    \psi(S_{T})=\left(K-\min_{1\leq i\leq d}S_{t}^{i}\right)_{+}.
    \]
\end{itemize}
The vector $(\lambda_{1},...,\lambda_{d})$ represents the weight of the
different assets within the basket. These weights may be negative to allow
to consider exchange options.  The variable $K$ denotes the strike price
and the vector $U$ corresponds to an upper barrier on the asset price at
maturity time $T$. Most of the time, the expectation $\E(\psi(S_{T}))$ needs
to be computed numerically by means of Monte Carlo methods.  In order to use
a systematic approach, the expectation $\E(\psi(S_{T}))$ is usually rewritten
as $\E(\phi(G))$ where $G$ is a standard normal random vector in
$\mathbb{R}^{d}$ and $\phi:\mathbb{R}^{d}\rightarrow\mathbb{R}^{+}$ is a
measurable and integrable function. To use our adaptive method described in
the next section, we need to transform the computation of $\E(\phi(G))$ into
the numerical computation of
\[
I(A)=\int_{[-A,A]^{d}}\phi(x)p(x)dx
\]
where $p(x)$ is the density of $G$ and $A\in\mathbb{R}_{+}$ is
chosen large enough to have a good approximation of $\E(\phi(G))$
by $I(A).$

\subsection{Key role of the delta}

When selling a financial derivative, the price of the option is obviously of
a primary interest but one should keep in mind that the original definition
of an option price goes back to the replicating theory. The price is actually
defined as the value at time $t=0$ of the replicating portfolio. Hence, the
price becomes fairly useless if we do not know how to implement the
replicating portfolio; hopefully, we precisely know how many stocks the
portfolio should carry on at any time $t$ and this quantity is given by the
famous \emph{Delta} of an option defined as the gradient of the option price
with respect to the spot vector
\[
\Delta=\nabla_{S_{0}}\exp(-rT)\E(\psi(S_{T})).
\]
There are several numerical methods for computing the delta of an option. The
most commonly used method is based on a finite difference approach because of
its automatic application. However, when coupled with a Monte Carlo method
it often yields a poor accuracy unless a very large number of samples is
used.

\section{Description of the adaptive algorithm}

\subsection{Quasi-Monte Carlo quadratures}

Any adaptive integration method relies on a quadrature rule designed for the
non-adaptive case. While Monte Carlo or quasi-Monte Carlo methods can deal
with integrands with no or little regularity, usual quadrature rules like
Gauss product rules are built for functions having a given  regularity or
belonging to particular bases.  These quadratures should not be too sensitive
to the dimensional effect and thus quadrature formulae based on interpolation
or approximation on different reduced size bases are considered. In the case
of $d-$dimensional Fourier bases on periodic smooth functions, Korobov spaces
\cite{key-KROMMER} are built using that the Fourier coefficients $a_{m}$
verify
\[
\left|a_{m}\right|\leq\frac{C}{(\widetilde{m_{1}}\widetilde{m_{2}}...\widetilde{m_{Q}})^{\beta}}
\]
where $\widetilde{m}=\max(1,\left|m\right|)$ and $\beta>1$ is linked
to the regularity of the integrand. The corresponding quadrature formulae
are lattice rules \cite{key-KROMMER,key-SLOAN} which are exact for the most
significant Fourier coefficients according to this decay. Non-periodic but smooth
functions can be periodized \cite{key-HELLUY} to still use these quadratures
but with an increasing constant of decay. 

In the case of polynomial approximations, Novak and Ritter \cite{key-NOVAK}
have obtained quadrature formulae, which are exact for polynomials with total
degree less than a given value. Based on the use of the control variates
method on piecewise interpolation polynomials, Atanassov and Dimov
\cite{key-ATANASSOV} built a numerical method reaching an optimal rate of
convergence for multivariate smooth functions with a fixed degree of
differentiation. 

The previous quadratures give highly accurate results for smooth functions
but they can yield a really poor accuracy for non-smooth functions. Monte
Carlo methods are somehow more stable because they are not sensitive to the
smoothness of the integrand.  The quadratures developed in \cite{key-MAIRE4}
combine the approximation on reduced Tchebychef polynomial basis and the use
of quasi-Monte Carlo points to build this approximation.  They are especially
efficient for very smooth functions but they can also handle pretty well only
continuous functions. They have been obtained after successive improvements
of an initial adaptive Monte Carlo method \cite{key-MAIRE1} via quasi-random
sequences \cite{key-MAIRE3} and the introduction of Tchebychef polynomial
basis of Korobov type \cite{key-MAIRE2}. We recall these formulae for a
multivariate function defined on the hypercube $[-1,1]^{d}$. For $m \in
{\mathbb N}$,
let $\widehat{m}=\max(1,m)$. We introduce the set 
\[
W_{d,q}=\left\{ m\in\mathbb{N}^{d}/\left(\widehat{m_{1}}...\widehat{m_{d}}\right)\leq q\right\} 
\]
which corresponds to the level $q$ of approximation. The reduced Tchebychef
polynomial approximation writes
\[
f(x_{1},x_{2},..,x_{d})\simeq\sum_{m\in W_{d,q}}b_{m}T_{m_{1}}(x_{1})T_{m_{2}}(x_{2})..T_{m_{d}}(x_{d})
\]
where the $L_{d,q}=card(W_{d,q})$ coefficients $b_{m}$ verify
\[
\left|b_{m}\right|\leq\frac{C_{1}}{(\widehat{m_{1}}\widehat{m_{2}}...\widehat{m_{d}})^{2L}}
\]
for a $C^{2L}$ function. To compute the numerical approximation of the
coefficients $b_{m},$ we fit the model 
\[
\sum_{m\in W_{d,q}}b_{m}T_{m_{1}}(x_{1})T_{m_{2}}(x_{2})..T_{m_{d}}(x_{d})
\]
to the observations of the function $f$ at some quadrature points
$P_{i}=(X_{i}^{(1)},...,X_{i}^{(d)})$ with $1\leq i\leq M.$ The choice of the
points $P_{i}$ is obviously crucial for the condition number $\kappa(A)$ of
the least-square matrix $A$.  We have proved in \cite{key-MAIRE4} that if
these points are independent random variables distributed according to the
multidimensional Tchebychef density
\[
w(x_{1},x_{2},..,x_{d})=\prod_{i=1}^{d}\frac{1}{\pi\sqrt{1-x_{i}^{2}}}1_{[-1,1]}(x_{i})
\]
then  $\kappa(A)$ goes to 1 when $M\rightarrow\infty$ at a Monte Carlo speed
$\frac{\sigma}{\sqrt{M}}$.  Moreover the variance $\sigma^{2}$ is  bounded by
one for any set $W_{d,q}$.  An even better choice to represent the density
$w$ is to use quasi-Monte Carlo points or quantization points
\cite{key-PAGES}.

We compute numerically the inverse of the matrix $A$ in order to obtain once
and for all quadrature formulae for each of the coefficients $b_{m}$ and
finally a quadrature formula for the integral itself. In the following, the
quadrature points are built using $\alpha\times L_{d,q}$ Halton points with
one additional point at each corner of the domain for a total of
$M=\alpha\times L_{d,q}+2^{d}$ points. In most situations, the parameter
$\alpha$ is chosen equal to $3$ which is sufficient to ensure a small value
for  $\kappa(A)$. The corner points are control points to detect a possible
change of regularity of the function $f$. Nevertheless, we will see in
Section 4, that more quadrature points may be necessary to detect this change
of regularity on the difficult example of digital options. The approximation
of
\[
I(f)=\int_{[-1,1]^{d}}f(x)dx
\]
is given by the quadrature formula
\[
Q_{d,\alpha,q}(f)=\sum_{i=1}^{\alpha\times
L_{d,q}+2^{d}}\omega_{i}f(X_{1}^{(i)},...,X_{d}^{(i)})
\]
which is obtained thanks to the corresponding approximations
$Q_{d,\alpha,q,b_{m}}(f)$ of the coefficients $b_{m}$ and the integration of
the approximation model.  We also denote  by $Q_{d,\alpha,q,R}(f)$ and
$Q_{d,\alpha,q,b_{m},R}(f)$ the relative quadrature formulae on a given
rectangle $R$ which will be used in the adaptive integration algorithm.

\subsection{Error indicators}

We keep the same error indicators than the one used in our previous paper \cite{key-CDL}. 
They rely on hierarchical quadratures for the integral of $f$ but also for some of the coefficients
of its weighted mean-square approximation. We select two sets $W_{d,q_{1}}$ and $W_{d,q_{2}}$ with
$1\leq q_{1}<q_{2}$, the corresponding quadrature formulae $Q_{d,\alpha,q_{1}}(f)$
and $Q_{d,\alpha,q_{2}}(f)$ for the integral of $f$ and also the quadrature
formulae for some of the leading coefficients in the approximation
model. These coefficients are the $d+1$ coefficients belonging to
the set $A_{d}$ of the coefficients $b_{m}$ for which all the indexes
$(m_{1},m_{2}..,m_{d})$ are equal to zero or only one is non-zero
and equal to one. Our error indicator $E_{d,\alpha,q_{1},q_{2},R}(f)$ is
\[
\left|Q_{d,\alpha,q_{1},R}(f)-Q_{d,\alpha,q_{2},R}(f)\right|+\sum_{m/b_{m}\in A_{d}}\left|Q_{d,\alpha,q_{1},b_{m},R}(f)-Q_{d,\alpha,q_{2},b_{m},R}(f)\right|
\]
on a given hyperrectangle $R.$ This indicator is more robust than the usual
indicator based on the comparison between only $Q_{d,\alpha,q_{1},R}(f)$ and
$Q_{d,\alpha,q_{2},R}(f)$. Indeed, it is less likely to happen that all the
estimators of the leading coefficients are close to each other but all wrong.
The approximate value of the integral on a given hyperrectangle is
$Q_{d,\alpha,q_{2},R}(f)$ and the approximate integral of $f$ is the sum of the
integrals over all the hyperrectangles.

\subsection{Splitting strategies}

Now, we describe the splitting strategies used in the numerical experiments.
At each step of the algorithm, the list of all the remaining hyper-rectangles
involved in the algorithm is stored and these hyperrectangles are sorted
according to their error indicator. Once the splitting is performed, the
hyperrectangle $R$ with the largest error indicator is removed from the list
and the hyperrectangles $R_{1}$ and $R_{2}$ are inserted in the list
according to their error indicators $E_{d,\alpha,q_{1},q_{2},R_{1}}(f)$ and
$E_{d,\alpha,q_{1},q_{2},R_{2}}(f)$.

\subsubsection{Fully adaptive splitting}

The first strategy named FAS is fully adaptive. It consists in trying the $d$
possible ways to divide the hyper-rectangle $R$ with largest error indicator
in two equal size pieces $R_{1}$ and $R_{2}$ along one of the axes and
keeping only the best splitting. The best splitting is the one for which 
\[
E_{d,\alpha,q_{1},q_{2},R_{1}}(f)+E_{d,\alpha,q_{1},q_{2},R_{2}}(f)
\]
is minimum among the $d$ possible choices.  Each step of this algorithm
requires $dN(\alpha L_{d,q}+2^{d})$ evaluations of the function $f$. The
computational cost of the method should also take into account the cost of
inserting the two newly created hyperrectangles into the list of all
hyperrectangles. This can be done efficiently by an insertion sort which has
a complexity of $O(N \log(N))$ for our algorithm with $N$ steps. In most
practicals examples, the evaluation of the function $f$ requires a lot of
atomic computations and hence it makes sense to neglect the cost the sort
compared to the number of function evaluations.

\subsubsection{Geometrical Random Splitting}

The FAS requires $d$ trials to find the optimal splitting. 
It is mainly interesting when the integrand has variations different
from several degrees of magnitude from one coordinate to another like for
instance for the two dimensional function $\cos(200x+y).$ When pricing vanilla
options, this kind of situation is unlikely to happen. Furthermore
and maybe more importantly, we have noticed in our previous work
\cite{key-CDL} that the convergence problems of our algorithm came from a too
fine splitting in one direction near the interface between the regions where
the regularity changes.

We propose another strategy named geometrical random splitting (GRS) in order
to reduce the computational costs and to solve the convergence problems by
finding more precisely this interface.  At each step of the algorithm, we
divide  the hyperrectangle with the largest error criterion  in two equal
size pieces $R_{1}$ and $R_{2}$ uniformly at random among the admissible
directions. The admissible directions are the ones having the larger length
among the axis of the hyperrectangle.  The initial hyperrectangle is always a
hypercube which means that we can never have a ratio more than two between
the different lengths of the axis of the  hyperrectangles occurring in the
mesh.  The cost of the algorithm is  $N(\alpha L_{d,q}+2^{d})$ and it is
stochastic which enables to compute some statistics on its results.

\section{Pricing Vanilla options in low dimensions}

In this section, we compare the FAS and the GRS for pricing vanilla
options on several examples already treated in \cite{key-CDL} and also on
examples from the general model. On all  numerical examples, the same
quadrature formulae will be used for the two methods. The number of
iterations of the algorithm will be $N=2000$ for the FAS and consequently $d$
times more for the GRS in order to keep the same complexity.

\subsection{Basket options}

Let us focus a little on basket options. The price of such options, with
payoffs only depending on the asset price at maturity time $T$, can be
expressed as an expectation $e^{-rT} \E(\psi(S_T))$. When the random vector
$S_T$ has a density $f_S$ with respect to the Lebesgue measure, computing this
expectation is easily turned into a numerical integration problem
\[
\E(\psi(S_T)) = \int_{[-\infty, \infty]^d} \psi(x) f_S(x) dx.
\]
When $d=1$, there is no need of a numerical integration as a closed formula
exists. In the Black--Scholes framework, we know that 
\[
S_T \stackrel{law}{=} \left(
S_{0}^{i}\exp\left\{(r-\frac{\sigma_{i}^{2}}{2})T+\sigma_{i}C_{i} \sqrt{T}
\; G\right\} \right)_{1 \le i \le d}
\]
where $G$ is a random normal vector with values in $\R^d$. Hence, 
\begin{align*}
  \E(\psi(S_T)) = \int_{[-\infty,\infty]^d}
  \psi\left(S_{0}^{i}\exp\left\{(r-\frac{\sigma_{i}^{2}}{2})T+\sigma_{i}C_{i}
  \sqrt{T} \; x\right\} , 1 \le i \le d \right) \frac{1}{(2\pi)^{d/2}}  e^{-|x|^2/2} dx.
\end{align*}
For the particular case of a call basket option, the payoff $\psi$ writes down
\[
\psi(s) = \left( \sum_{i=1}^d \lambda_i s_i - K \right)_+.
\]
Therefore, if we denote by $V(T,K)$ the price of the call basket option, we can write
\[
V(T,K) = e^{-rT} \frac{1}{(2\pi)^{d/2}} \int_{[-\infty,\infty]^d} \left(
\sum_{i=1}^d \lambda_i
S_{0}^{i}\exp\left\{(r-\frac{\sigma_{i}^{2}}{2})T+\sigma_{i}C_{i} \sqrt{T}
\; x\right\} - K\right)_+  e^{-|x|^2/2} dx.
\]
Using the equality $(s - K)_+ - (K - s)_+ = s - K$ and letting the price of
the corresponding put basket option be $U(T,K)=e^{-r T} \E\left((K -
\sum_{i=1}^d \lambda_i S_T^i)_+\right)$, we obtain the so--called
call put parity relationship 
\[
V(T,K)-U(T,K)=\sum_{i=1}^d \lambda_i -K\exp(-rT).
\]

In general, this formula is used to obtain the hardest price to compute (in
terms of variance) between the call or the put option prices from the other.
Here, we independently compute the call and put option prices and use
this formula as a criterion of accuracy. The infinite integral is truncated
and we let
\[
V(T,K,A)=e^{-rT} \frac{1}{(2\pi)^{d/2}} \int_{[-A,A]^d} \left(
\sum_{i=1}^d\lambda_i
S_{0}^{i}\exp\left\{(r-\frac{\sigma_{i}^{2}}{2})T+\sigma_{i}C_{i} \sqrt{T} \;
x\right\} - K\right)_+  e^{-|x|^2/2} dx
\]
be the truncated estimation of $V(T,K)$. The truncated approximation of the
put basket option price  $U(T,K,A)$ is defined in a similar way. Note that
these multi--dimensional integrals are truncated on a square domain.
In all the following examples dealing with basket options, we fix the weights
of the basket as $\lambda_i = \frac{1}{d}$ so that the baskets are all
homogeneous and their weights sum up to one.

\subsubsection{Examples in dimension 2}

The numerical tests considered in this paragraph are very similar to the ones
treated in \cite{key-DEELSTRA}. Three different values for the strike price
are tested, one at the money $K_{1}=100,$ one out $K_{2}=127.80$ and one
completely out $K_{3}=300.$ These examples were already studied in
\cite{key-CDL} and the conclusion was that the FAS was outperforming Monte
Carlo and quasi-Monte Carlo integration. Hence, we only compare the GRS to
the FAS. Moreover, the new algorithm is only run once as we have a
deterministic error criterion based on the call-put parity formula.

The GRS is used with quadrature formulae of degrees $18$ and $24$. The latter
corresponds to a number of function evaluations of
$2000\times403\times2\simeq1.6\times10^{6}.$ We give in Tables
\ref{tab_2DBasket_ex1} and \ref{tab_2DBasket_ex2}, the values of the
truncated estimations for $A_{1}=12$ and $A_{2}=13$. We also compute 
\[
C(T,K,A)=\left|V(T,K,A)-U(T,K,A)-\frac{1}{2}S_0^{(1)}-\frac{1}{2}S_0^{(2)}+K\exp(-rT)\right|
\]
for these reference values and denote by $C_{old}(T,K,A)$ the same error
indicator using the FAS.  Truncating the domain of a real valued standard
normal random variable to $[-12,12]$ may seem far too large. This is true
from a Monte-Carlo point of view and we have actually also run some tests
with more conventional values such as $A=5$. These tests were already showing
a far better accuracy than a crude Monte Carlo; the results obtained on the
examples of Table~\ref{tab_2DBasket_ex1} already had $4$ accurate digits.
However, since we are targeting to compute sensibilities, we need a better
accuracy, which explains the choices of the parameter $A$.

\begin{table}[ht]
\centering\begin{tabular}{|c|c|c|c|c|}
  \hline 
  $\ $ & $V$ & $U$ & $C$ & $C_{old}$\tabularnewline
  \hline 
  $(K_{1},A_{1})$ & 28.49407706 & 14.564874729 & $2\times10^{-8}$  & $1\times10^{-8}$\tabularnewline
  \hline 
  $(K_{1},A_{2})$ & 28.49407708 & 14.564874726 & $2\times10^{-10}$  & $1\times10^{-8}$\tabularnewline
  \hline 
  $(K_{2},A_{1})$ & 18.85549194 & 28.853971355 & $2\times10^{-8}$ & $2\times10^{-8}$\tabularnewline
  \hline 
  $(K_{2},A_{2})$ & 18.85549196 & 28.853971353 & $2\times10^{-9}$ & $7\times10^{-9}$\tabularnewline
  \hline 
  $(K_{3},A_{1})$ & 1.810536572 & 160.02292952 & $2\times10^{-8}$ & $2\times10^{-8}$\tabularnewline
  \hline 
  $(K_{3},A_{2})$ & 1.810536593 & 160.02292952 & $4\times10^{-10}$ & $4\times10^{-10}$\tabularnewline
  \hline 
\end{tabular}
\caption{Basket option with parameters $T=3$, $r=0.05$, $d=2$ \newline$S_0^{(1)}=S_0^{(2)}=50$, $\sigma_1=\sigma_2=0.4,\rho=0.3$}\label{tab_2DBasket_ex1}
\end{table}
On this first example, we obtain a very good accuracy of at least
$8$ digits on all the computations of the prices of call and put
options. These values are both validated by the parity
call put formula and by the comparison between the truncated approximations
for the two different values of $A.$ There is no significant difference between
the two splitting strategies. 
\begin{table}[ht]
\centering\begin{tabular}{|c|c|c|c|c|}
  \hline 
  $\ $ & $V$ & $U$ & $C$ & $C_{old}$\tabularnewline
  \hline 
  $(K_{1},A_{1})$ & 20.04091112 & 6.1117087694 & $2\times10^{-9}$  & $1\times10^{-8}$\tabularnewline
  \hline 
  $(K_{1},A_{2})$ & 20.04091112 & 6.1117087676 & $1\times10^{-10}$  & $1\times10^{-8}$\tabularnewline
  \hline 
  $(K_{2},A_{1})$ & 8.915343209 & 18.913822596 & $7\times10^{-10}$  & $2\times10^{-8}$\tabularnewline
  \hline 
  $(K_{2},A_{2})$ & 8.915343211 & 18.913822598 & $5\times10^{-10}$ & $7\times10^{-9}$\tabularnewline
  \hline 
  $(K_{3},A_{1})$ & 0.021755879 & 158.23414880 & $4\times10^{-10}$ & $2\times10^{-8}$\tabularnewline
  \hline 
  $(K_{3},A_{2})$ & 0.021755880 & 158.23414880 & $6\times10^{-11}$ & $4\times10^{-10}$\tabularnewline
  \hline 
\end{tabular}
\caption{Basket option with parameters $T=3$, $r=0.05$, $d=2$\newline$S_0^{(1)}=S_0^{(2)}=50$, $\sigma_1=\sigma_2=0.2,\rho=0.7$}\label{tab_2DBasket_ex2}
\end{table}
The two versions of the algorithm are still very efficient on an
example with a smaller volatility and more correlated assets. The accuracy
is even better than on the first example and we note that the GRS is now one or two digits more accurate than the FAS.
\begin{figure}[h]
  \begin{center}
    \subfigure[Geometrical Random Splitting]
    {\includegraphics[width=0.495\textwidth]{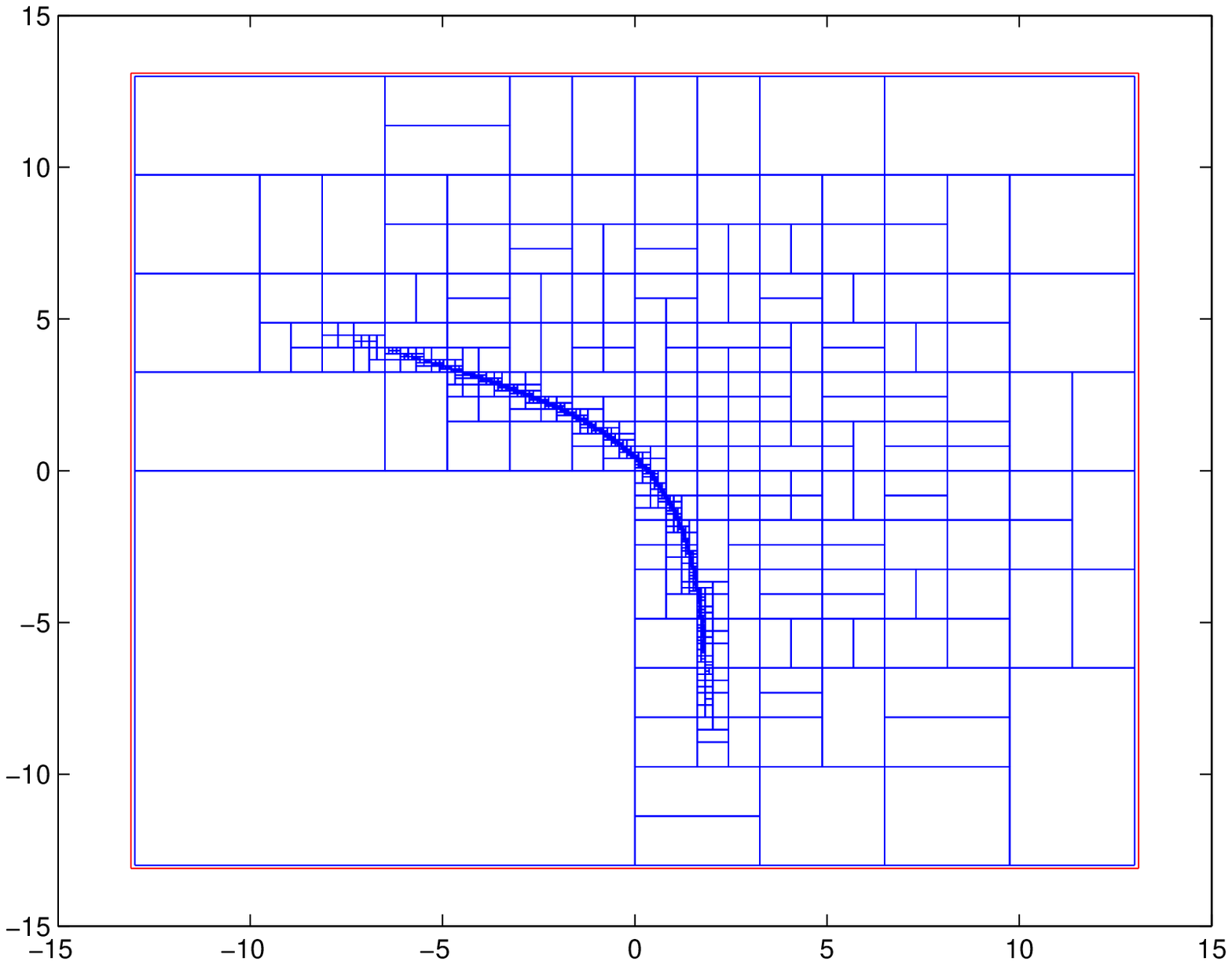}}
    \subfigure[Fully Adaptive Splitting]
    {\includegraphics[width=0.495\textwidth]{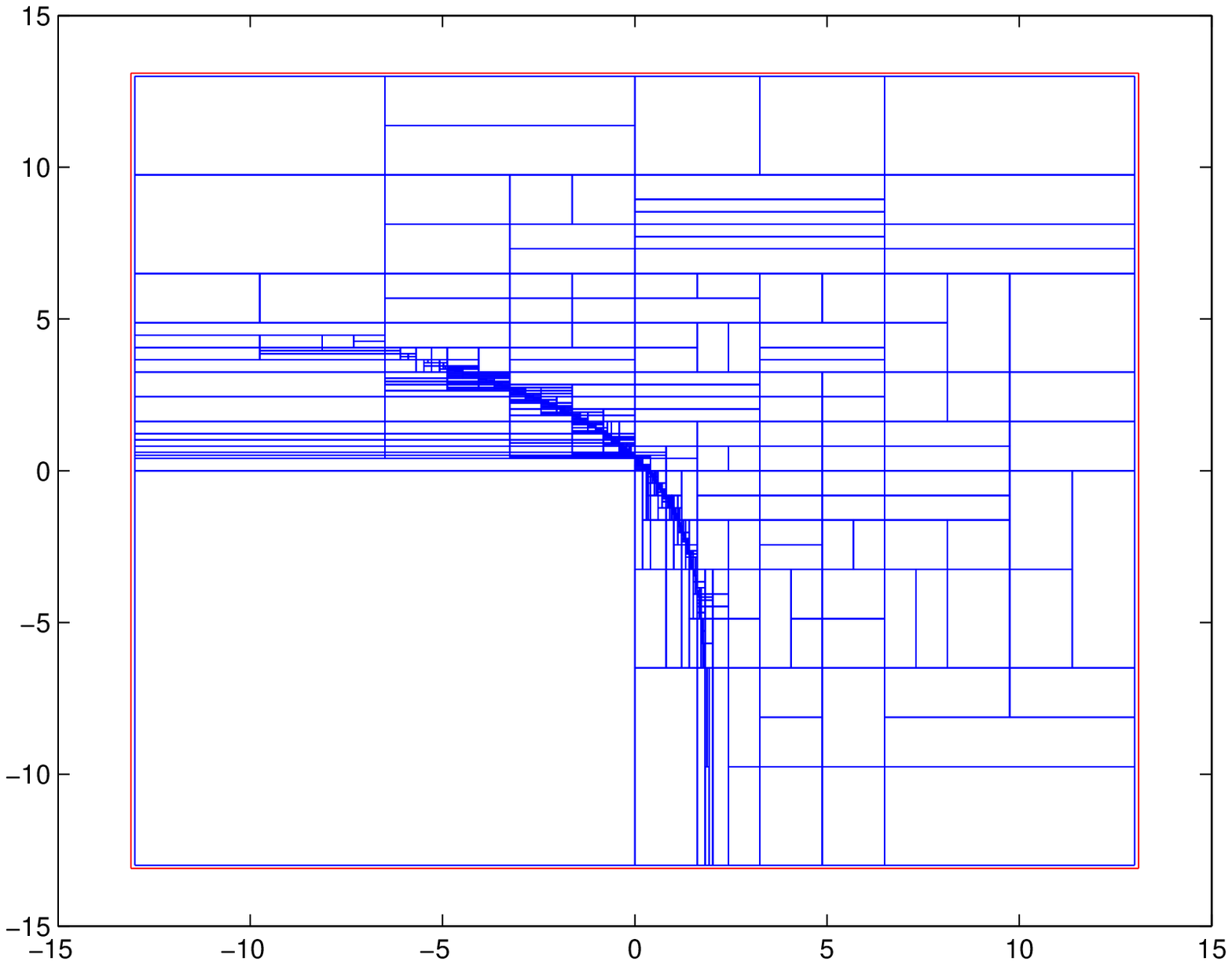}}
    \caption{Mesh for option of Table~\ref{tab_2DBasket_ex1} with $K = K_1 = 100$}
    \label{fig_2DBasket_ex1}
  \end{center}
\end{figure}
It is also interesting to compare the meshes obtained for the call options
for the two splitting strategies.  In Figure \ref{fig_2DBasket_ex1}, we plot these meshes for the first example with parameters $(K,A)=(K_{1},A_{2})$.  
We observe that the refinement is done mainly near the interface separating the region where the
function vanishes and the region where it is positive. 
On these figures, we can see that the GRS is more robust and
more able to follow the interface than the FAS.

\subsubsection{Examples in dimension $3$ and $4$}

In this paragraph, we consider two examples of basket options on independent
assets: one example in dimension $3$ and one in dimension $4$.
We will run each example with two different strike prices.

We give in Table \ref{tab_3DBasket} the values of the truncated
estimations for $A_{1}=12$ and $A_{2}=13$ for a $3$ dimensional basket
option. The indicator error for this example is based on $E_{3,3,18,24,R}(f)$
and the number of function evaluations is given by $2000\times1592\times3\simeq9.5\times10^{6}.$

\begin{table}[ht]
\centering\begin{tabular}{|c|c|c|c|c|}
  \hline 
  $\ $ & $V$ & $U$ & $C$ & $C_{old}$\tabularnewline
  \hline 
  $(K_{1},A_{1})$ & 14.80805242 & 2.2717704262 & $1\times10^{-7}$  & $4\times10^{-6}$\tabularnewline
  \hline 
  $(K_{1},A_{2})$ & 14.80805257 & 2.2717705311 & $7\times10^{-8}$   & $7\times10^{-7}$\tabularnewline
  \hline 
  $(K_{2},A_{1})$ & 2.927052540 & 16.212009773 & $6\times10^{-8}$  & $8\times10^{-5}$\tabularnewline
  \hline 
  $(K_{2},A_{2})$ & 2.927053375 & 16.212010568 & $2\times10^{-8}$  & $5\times10^{-5}$\tabularnewline
  \hline
\end{tabular}
 \caption{$d=3,T=3,r=0.05,\rho=0.3,S_0^{(1)}=S_0^{(2)}=S_0^{(3)}=30\newline\sigma_{1}=\sigma_{2}=\sigma_{3}=0.2
 $, $K_1=90, K_2=120$}\label{tab_3DBasket}
\end{table}

We observe that the GRS is a lot more accurate than the FAS.
The accuracy is about $7$ digits on the at the money option and $8$ on the
one out of the money. On the out of the money option, the FAS hardly
achieves $5$ digits of accuracy. \\

We give in Table \ref{tab_4DBasket} the values of the truncated estimations
for $A_{1}=5$ and $A_{2}=6$. The indicator error for this example is based on
$E_{4,3,18,24,R}(f)$ and the number of function evaluations is given by
$2000\times9121\times4\simeq4.4\times10^{7}$.

\begin{table}[ht]
\centering\begin{tabular}{|c|c|c|c|c|c|}
  \hline 
  $\ $ & $V$ & $U$ & $C$ &  $C_{old}$\tabularnewline
  \hline 
  $(K_{1},A_{1})$ & 4.22830628 & 0.32667437 & $1\times10^{-5}$  & $3\times10^{-5}$\tabularnewline
  \hline 
  $(K_{1},A_{2})$ & 4.22832492 & 0.32667871 & $1\times10^{-7}$  & $4\times10^{-4}$\tabularnewline
  \hline 
  $(K_{2},A_{1})$ & 0.16841321 & 5.77905377 & $7\times10^{-6}$ & $7\times10^{-4}$\tabularnewline
  \hline 
  $(K_{2},A_{2})$ & 0.16842047 & 5.77906874 & $6\times10^{-8}$ & $6\times10^{-4}$ \tabularnewline
  \hline
\end{tabular}
 \caption{$d=4,T=1,r=0.05,\rho=0.3,S_0^{(1)}=S_0^{(2)}=S_0^{(3)}=S_0^{(4)}=20\newline\sigma_{1}=\sigma_{2}=\sigma_{3}=\sigma_{4}=0.1$, $K_1=80, K_2=90$}\label{tab_4DBasket}
\end{table}

On this last example, the improvement of the algorithm is even more
impressive. The accuracy is about $7$ digits with the GRS compared to
hardly $4$ with the FAS. The difference between the two values of  the
criterion $C$ for the two values of $A$ is not due to a wrong computation of
the integrals but simply because the integration domain was  truncated too
much when $A=A_{1}.$

\subsection{Put on minimum options} 

For the pricing of basket options, our approximation criterion was relying on
the parity call put formula which does not hold for on minimum or digital
options. It required also the pricing of a call and a put to be computed. To
define error criteria for general options, we exploit the stochastic nature
of the GRS.  Instead of running only once this algorithm, we run it ten times
only on the call pricing and compute its average $\bar{V}$,  its empirical
variance $s^2_v$ and our error criterion will be simply $Err=s_v$.  We also
give the approximation  $V_{old}$ obtained with the FAS and the estimation
$V_{MC}$ obtained using the crude Monte Carlo method for the computation of
the price as the expectation of a normal random vector as described in
Section~\ref{sec_model}.  This Monte Carlo estimation is accurate up to three
digits in relative error and we also give the number of samples required to
reach such an accuracy.  \subsubsection{Examples in dimension 2}

Numerical results are given in Table \ref{tab_2DWorst} for both cases $A_{1}=12$ and $A_{2}=15.$ 
\begin{table}[H]
\centering\begin{tabular}{|c|c|c|c|c|}
  \hline 
  $\ $ & $\bar{V}$ & $V_{old}$ & $V_{MC}$ & $Err$ \tabularnewline
  \hline 
  Ex1, $A_{1}$ & 2.10306340730 & 2.10306339974 & 2.104291  & $1.5\times10^{-10}$\tabularnewline
  \hline 
  Ex1, $A_{2}$ & 2.10306340508 & 2.10306346643 & 2.104291  & $1.2\times10^{-10}$\tabularnewline
  \hline 
  Ex2, $A_{1}$ & 6.32237986596 & 6.32237987060 & 6.325378 & $2.2\times10^{-10}$\tabularnewline
  \hline 
  Ex2, $A_{2}$ & 6.32237986541 & 6.32237985738 & 6.325378 & $1.2\times10^{-10}$\tabularnewline
  \hline
\end{tabular}
\caption{Put on minimum options in dimension $d=2$ \newline
Ex1:
$T=1,r=0.05,S_0^{(1)}=S_0^{(2)}=50,\rho=0.1,K=45,\sigma_{1}=\sigma_{2}=0.2$\newline
Ex2:
$T=1,r=0.05,S_0^{(1)}=S_0^{(2)}=50,\rho=0.9,K=55,\sigma_{1}=\sigma_{2}=0.2$}\label{tab_2DWorst}
\end{table}
We observe that the GRS provides an
estimation of the price which is accurate up to 10 digits. Indeed the
integrals for the two values of $A$ are computed with an accuracy of more
than 10 digits and they have 10 common digits. The estimations obtained with
the FAS are slightly less accurate up to 8 digits. 
\begin{figure}[h]
  \begin{center}
    \subfigure[Geometrical Random Splitting]
    {\includegraphics[width=0.495\textwidth]{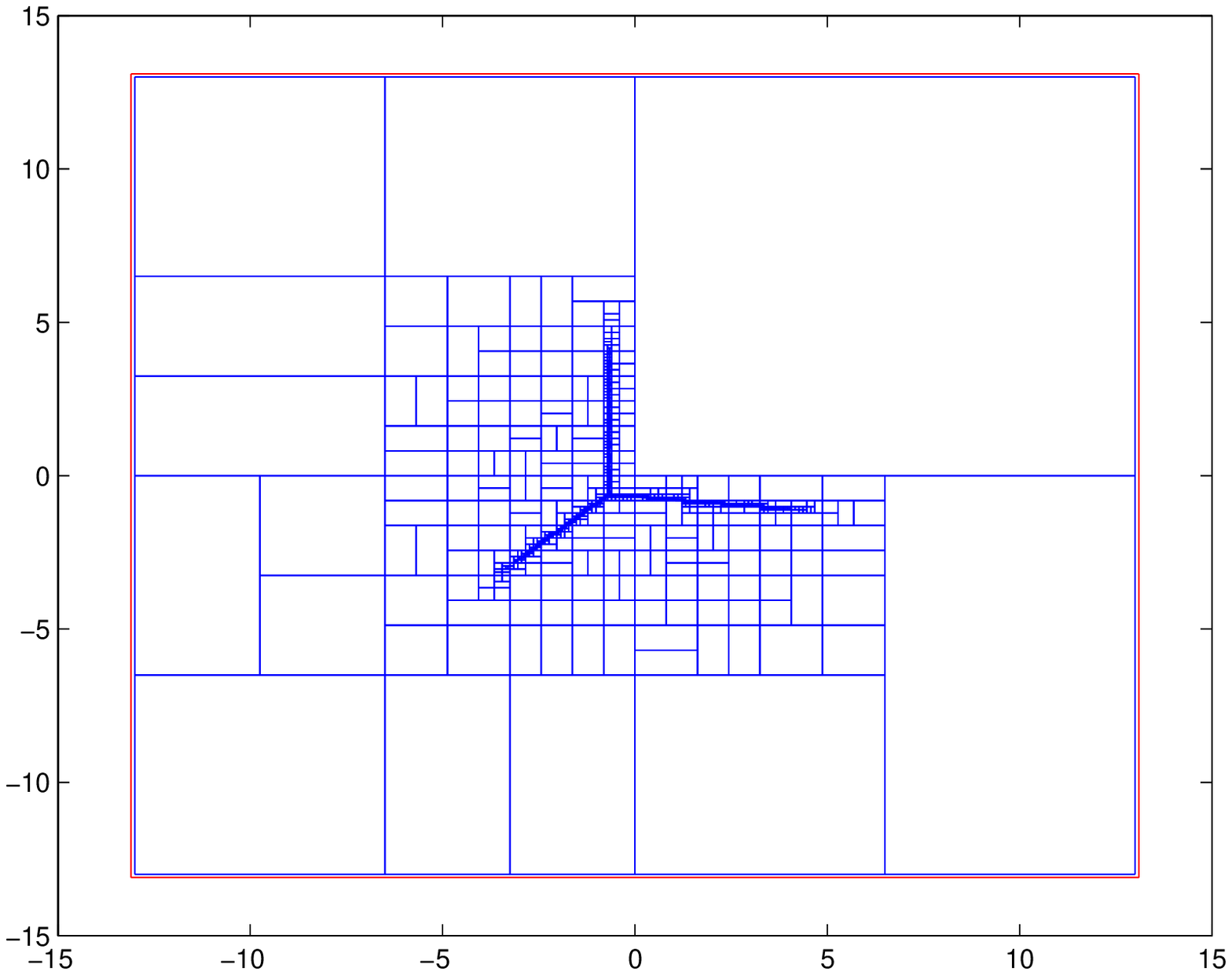}}
    \subfigure[Fully Adaptive Splitting]
    {\includegraphics[width=0.495\textwidth]{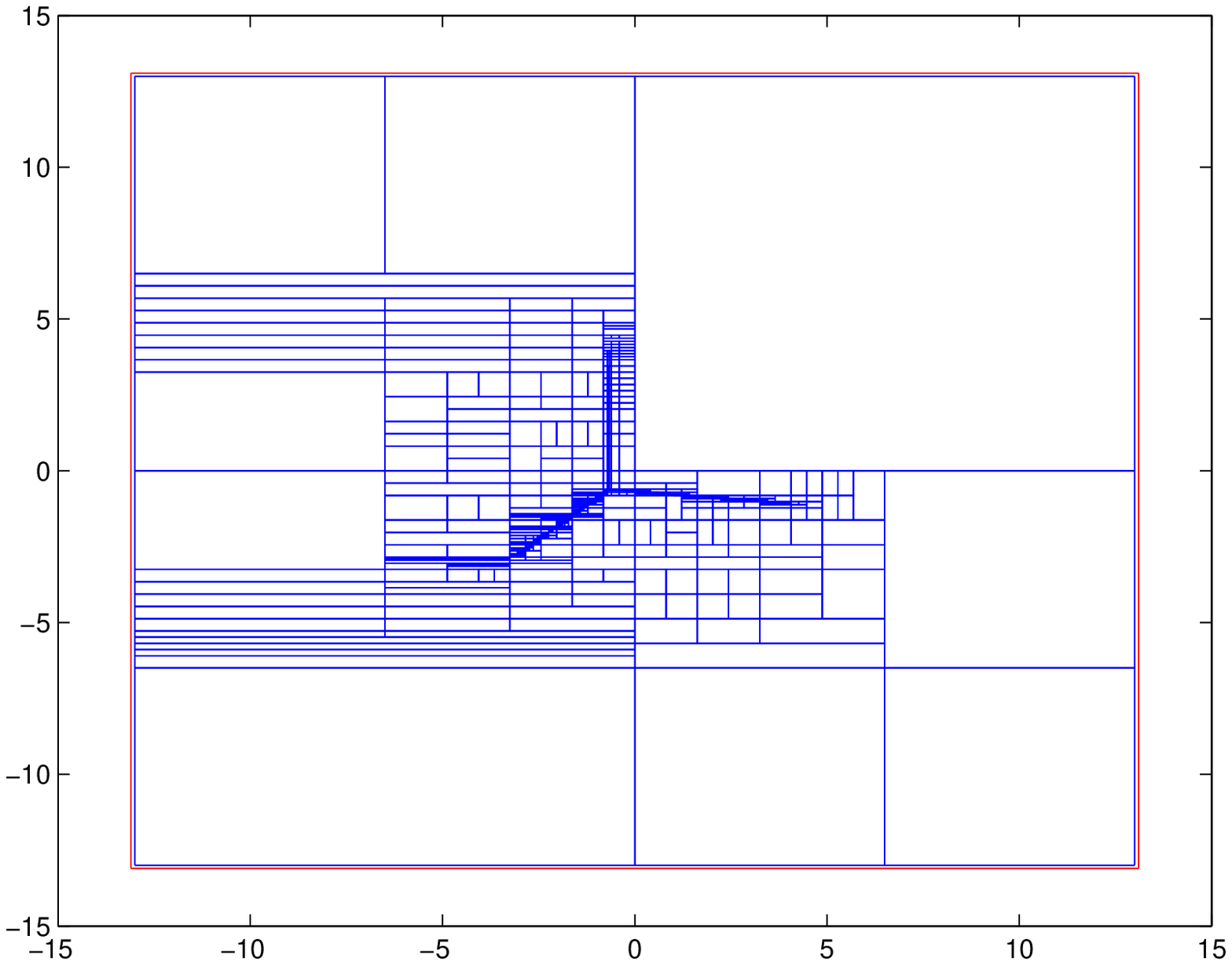}}
    \caption{Mesh for option of Table~\ref{tab_2DWorst} with $K=45$ (Ex1)}
    \label{fig_2DWorst}
  \end{center}
\end{figure}
In Figure \ref{fig_2DWorst}, we plot the two meshes for Ex1 in the case
$A=A_{1}.$ We observe on the GRS that the refinement is still done near the
interface separating the region where the function vanishes and the region
where it is positive but also near a line in the positive part. This
additional line is explained by a lack of regularity of the integrand due to
the minimum function in its definition. On these two meshes, we can see that
once again the GRS is better than the FAS to follow the interfaces where the regularity
changes. For both examples, one run of our algorithm requires
$1.6\times10^{6}$ function evaluations meanwhile there are $4.6\times10^{7}$
and $1.6\times10^{8}$ Monte Carlo samples for Ex1 and Ex2, respectively.  Even
if we take into account that 10 runs are used in order to calculate the error
indicator, the number of function evaluations for both methods (GRS and crude
Monte Carlo) are quite similar. This remains true in the following examples
in higher dimensions.

\subsubsection{Examples in dimensions 3 and 4}
We give in Table \ref{tab_3DWorst} the values of the truncated estimations
for $A_{1}=12$ and $A_{2}=15$ in dimension 3. 
\begin{table}[H]
\centering\begin{tabular}{|c|c|c|c|c|}
  \hline 
  $\ $ & $\bar{V}$ & $V_{old}$ & $V_{MC}$ & $Err$ \tabularnewline
  \hline 
  Ex3, $A_{1}$ & 2.89538461 & 3.021407 & 2.898180  & $6.3\times10^{-8}$\tabularnewline
  \hline 
  Ex3, $A_{2}$ & 2.89538389 & 2.909697 & 2.898180  & $3.1\times10^{-8}$\tabularnewline
  \hline 
  Ex4, $A_{1}$ & 6.85473710 & 6.802445 & 6.854480 & $6.3\times10^{-8}$\tabularnewline
  \hline 
  Ex4, $A_{2}$ & 6.85473692 & 6.795585 & 6.854480 & $6.3\times10^{-8}$\tabularnewline
  \hline 
\end{tabular}
\caption{Results for put on minimum options in dimension $d=3$\newline with
$r=0.05$, $S_0^{(1)}=S_0^{(2)}=S_0^{(3)}=50$, $T=1$ and
$\sigma_{1}=\sigma_{2}=\sigma_{3}=0.2$. \newline
Ex3: $\rho=0.1,K=45$\newline
Ex4: $\rho=0.9,K=55$}\label{tab_3DWorst}
\end{table}
We can see that the GRS is still very efficient and achieves an accuracy of
about 7 digits. Now, the FAS gives poorly accurate results and is also
clearly less efficient than crude Monte Carlo integration. This confirms the
difficulties for the FAS to capture the interface when the dimension of the
problem increases.  We give in Table \ref{tab_4DWorst}, the values of the
truncated estimations for $A_{1}=12$ and $A_{2}=15$ in dimension 4. The
conclusions are the same than in dimension 3. The FAS fails to converge while
the GRS is still very accurate.  On all examples, whatever the dimension is
(up to 4), this latter strategy outperforms the crude Monte Carlo method as
we have at least 6 digits of accuracy instead of 3 for a similar complexity.
\begin{table}[H]
\centering\begin{tabular}{|c|c|c|c|c|}
  \hline 
  $\ $ & $\bar{V}$ & $V_{old}$ & $V_{MC}$ & $Err$ \tabularnewline
  \hline 
  Ex5, $A_{1}$ & 3.567971 & 3.560892 & 3.574086  & $6.3\times10^{-7}$\tabularnewline
  \hline 
  Ex5, $A_{2}$ & 3.567971 & 3.322140 & 3.574086  & $3.1\times10^{-7}$\tabularnewline
  \hline 
  Ex6, $A_{1}$ & 7.212993 & 6.693007 & 7.215822 & $3.1\times10^{-7}$\tabularnewline
  \hline 
  Ex6, $A_{2}$ & 7.212994 & 7.305357 & 7.215822 & $3.1\times10^{-7}$\tabularnewline
  \hline
\end{tabular}
\caption{Results on minimum options in dimension $d=4$\newline with $r=0.05$,
$S_0^{(1)}=S_0^{(2)}=S_0^{(3)}=S_0^{(4)}=50$, $T=1$ and
$\sigma_{1}=\sigma_{2}=\sigma_{3}=\sigma_{4}=0.2$ \newline
Ex5: $\rho=0.1,K=45$ \newline
Ex6: $\rho=0.9,K=55$ }\label{tab_4DWorst}
\end{table}

\subsection{Digital options} As we will see in the following numerical
examples, the pricing of digital options is somehow a much harder problem
than the two previous ones. In fact, it is close to the problem of looking
for a small subdomain (where the function is positive) somewhere in a large
hyperrectangle.  This is well illustrated in Figure \ref{fig_2DDigit_B} and
even more in Figure \ref{fig_2DDigit_C} where this small subdomain is
delimited by three lines.  Especially during the first iterations of the
algorithm, it may happen that a hyperrectangle of our mesh contains a part of
this subdomain but with no quadrature points lying in it. In this case, our
error indicator is zero because the function is considered as the null
function and consequently this hyperrectangle is never split again. For
basket or put on minimum options, we had been able to handle this sampling
problem by putting additional quadrature points in each corner of the
hyperrectangle. This trick is inefficient in the case of digital options
because the subdomain with positive values is small and has a completely
different form. In order to solve at least partially these convergence
problems, we have increased the parameter $\alpha$ in the number of
quadrature points $M=\alpha\times L_{d,q}+2^{d}$ to obtain a larger
probability to detect a non-zero value of the function inside the
hyperrectangles. Finally, we also compute the median $Med(V)$ of 10 runs
of the algorithm because it is a more robust estimator in case of false
convergence.

\subsubsection{Examples in dimension 2}
Our numerical results 
for both cases $A_{1}=12$ and $A_{2}=15$ are given in Table \ref{tab_2DDigital}.

\begin{table}[H]
\centering\begin{tabular}{|c|c|c|c|c|}
  \hline 
  $\ $ & $\bar{V}$ & $Med(V)$ & $V_{MC}$ & $Err$ \tabularnewline
  \hline 
  Ex7, $A_{1}$, $\alpha=3$ & 2.30072052 &  2.30072041 & 2.299709 & $3.1\times10^{-7}$\tabularnewline
  \hline 
  Ex7, $A_{2}$, $\alpha=3$ & 2.30071826 & 2.30071825 & 2.299709  & $3.1\times10^{-7}$\tabularnewline
  \hline 
  Ex8, $A_{1}$, $\alpha=3$ & 0.12540527 &  0.15675651 & 0.15600  & $6.3\times10^{-2}$\tabularnewline
  \hline 
  Ex8, $A_{2}$, $\alpha=3$ & 0.13848118 & 0.15675549 & 0.15600  & $2.8\times10^{-2}$\tabularnewline
  \hline 
  Ex8, $A_{1}$, $\alpha=15$ & 0.15693827 & 0.15693825 & 0.15600 & $3.1\times10^{-8}$\tabularnewline
  \hline 
  Ex8, $A_{2}$, $\alpha=15$ & 0.15681002 & 0.15675531 & 0.15600 & $9.5\times10^{-5}$\tabularnewline
  \hline
\end{tabular}
\caption{Results for digital call options in dimension $d=2$\newline with
$T=1,r=0.05,S_0^{(1)}=S_0^{(2)}=50,U_1=U_2=60,\sigma_{1}=\sigma_{2}=0.2$ \newline
Ex7: $\rho=0.1,K=45$ \newline
Ex8: $\rho=0.9,K=55$}\label{tab_2DDigital}
\end{table}

\begin{figure}[H]
  \begin{center}
    \subfigure[Geometrical Random Splitting, $A_1$]
    {\includegraphics[width=0.55\textwidth]{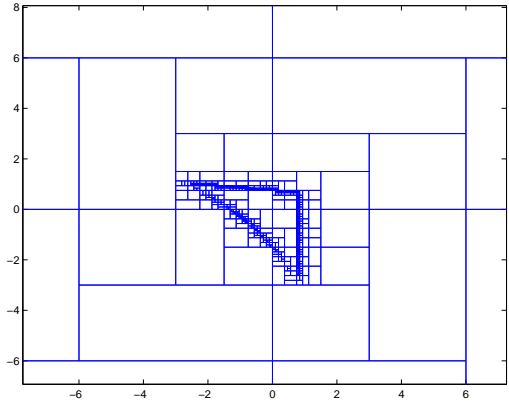}}
    \caption{Mesh for option Ex7 of Table~\ref{tab_2DDigital} $\alpha=3$}
    \label{fig_2DDigit_A}
  \end{center}
\end{figure}

On Ex7, we can see that we have no problem of convergence even with $\alpha=3$.
The price is approximately equal to $2.300718$, mean and median are at least
6 digits close in all cases which is  confirmed by the error indicator. We
also have 3 common digits with the Monte Carlo estimator. The mesh obtained
in Figure \ref{fig_2DDigit_A} shows that the subdomain where the function is
positive is roughly a small triangle with a surface $\xi$ lower than 8. The
ratio $A_{2}^2/\xi\simeq28$ is still small enough to avoid sampling problems. 

\begin{figure}[H]
  \begin{center}
    \subfigure[Geometrical Random Splitting, $A_2$]
    {\includegraphics[width=0.495\textwidth]{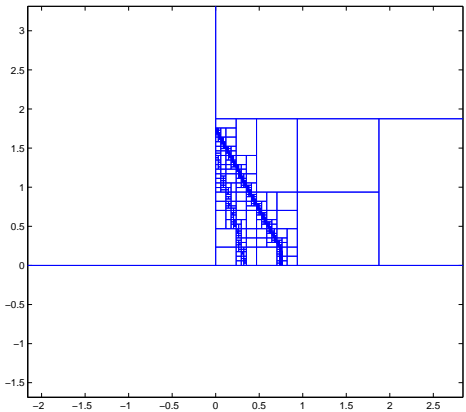}}
    \subfigure[Geometrical Random Splitting, $A_1$]
    {\includegraphics[width=0.495\textwidth]{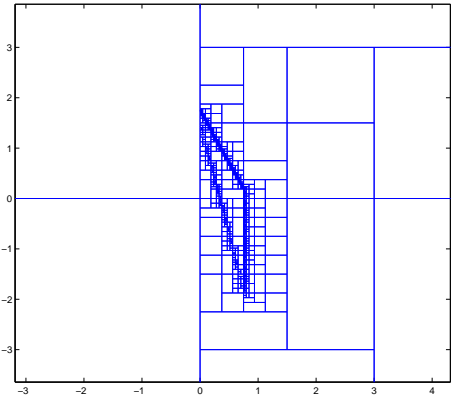}}

    \caption{Mesh for option Ex8 of Table~\ref{tab_2DDigital} $\alpha=3$}
    \label{fig_2DDigit_B}
  \end{center}
\end{figure}

\begin{figure}[H]
  \begin{center}
    \subfigure[Geometrical Random Splitting, $A_1$]
    {\includegraphics[width=0.55\textwidth]{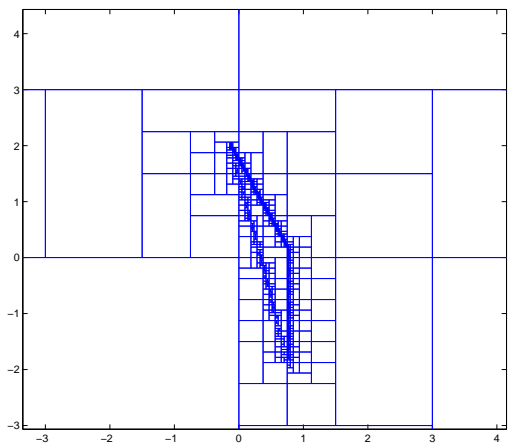}}
    \caption{Mesh for option Ex8 of Table~\ref{tab_2DDigital} $\alpha=15$}
    \label{fig_2DDigit_C}
  \end{center}
\end{figure}

The situation is completely different for Ex8 when  $\alpha=3$.  Mean and
median have no more than 2 common digits, the error indicator is also close
to $0.01$. However the median has 3 common digits with the crude Monte Carlo
estimator which means that the big convergence problems do not happen so
often. We plot in Figure \ref{fig_2DDigit_B} two examples of meshes
obtained with $A_{1}$ and $A_{2}$ in situations where the algorithm did not
converge. We observe  that the refinement is incomplete in both cases and
that  the one obtained with $A_{2}$ misses half of the region of interest.
If we add quadrature points by taking $\alpha=15$, we recover a high accuracy
of 7 digits on the integrals at least when  $A_{1}=12$. The subdomain we
obtain in this case in Figure \ref{fig_2DDigit_C} has now a triangular
shape. Its surface is now roughly equal to 1 which explains the convergence
problems we encountered when $\alpha=3$ as the ratio $A_{2}^2/\xi\simeq225.$  

For $\alpha=15$, one run of the GRS requires $8\times10^{6}$ function
evaluations, which is 5 times more than in the case $\alpha=3$ meanwhile there
are $4.4\times10^{7}$ and $1.4\times10^{6}$ Monte Carlo samples for
respectively Ex7 and Ex8.

\subsubsection{Examples in dimension 3}

We give in Table \ref{tab_3DDigital}, the values of the truncated
estimations for $A_{2}=12$ and $A_{2}=15.$ 
\begin{table}[H]
\centering\begin{tabular}{|l|c|c|c|c|}
  \hline 
  $\ $ & $\bar{V}$ & $Med(V)$ & $V_{MC}$ & $Err$ \tabularnewline
  \hline 
  Ex9,  $\ A_{1}$, $\alpha=20$  & 1.64950182 &1.64950187 & 1.646800  & $2.5\times10^{-6}$\tabularnewline
  \hline 
  Ex9,  $\ A_{2}$, $\alpha=20$  & 1.64948232 &1.64948236 & 1.646800  & $1.9\times10^{-6}$\tabularnewline
  \hline 
  Ex10, $A_{1}$, $\alpha=40$  & 0.09316076 &0.09316072 & 0.093269  & $3.1\times10^{-7}$\tabularnewline
  \hline 
  Ex10, $A_{2}$, $\alpha=40$  & 0.09307638 &0.09316133 & 0.093269 & $2.5\times10^{-4}$\tabularnewline
  \hline
\end{tabular}
\caption{Results on digital options in dimension $d=3$\newline with
$r=0.05,T=1,S_0^{(1)}=S_0^{(2)}=S_0^{(3)}=50,U_1=U_2=U_3=60$ \newline
Ex9: $\rho=0.1,K=45$\newline
Ex10: $\rho=0.9,K=55$}\label{tab_3DDigital}
\end{table}

The number of function evaluations are $6.3\times10^{7}$ for $\alpha=20$ and
$1.3\times10^{8}$ for $\alpha=40$ meanwhile there are $3.0\times10^{7}$ and
$7.3\times10^{6}$ Monte Carlo runs respectively for Ex9 and Ex10.  On Ex9, we
need to increase $\alpha$ to converge up to 5 or 6 digits but the number of
function evaluations is still comparable to the size of the Monte Carlo
sample. On Ex10 with $A_2$, the Monte Carlo method has a smaller variance and
the GRS a weaker accuracy which makes these two methods comparable in terms
of efficiency.  However, with $A_1$, the results obtained with the GRS seems
very accurate in comparison with those obtained with the Monte Carlo Method.
A more detailed discussion about the influence of the initial volume of the
hypperrectangle will be given in dimension 4.

\subsubsection{Examples in dimension 4}

We give in Table \ref{tab_4DDigital}, the values of the truncated estimations
for  the usual values $A_{1}=12$, $A_{2}=15$ and also for two smaller values
$A_{3}=5$ and $A_{4}=6.$
\begin{table}[H]
\centering\begin{tabular}{|c|c|c|c|c|}
  \hline 
  $\ $ & $\bar{V}$ & $Med(V)$ & $V_{MC}$ & $Err$ \tabularnewline
  \hline 
  Ex11, $A_{1}$, $\alpha=30$ & 1.22934667 & 1.22934613 & 1.228935  & $3.1\times10^{-6}$\tabularnewline
  \hline 
  Ex11, $A_{2}$, $\alpha=30$ & 1.22934272 & 1.22934341 & 1.228935  & $3.1\times10^{-6}$\tabularnewline
  \hline 
  Ex12, $A_{1}$, $\alpha=40$ & 0.04827239 & 0.04826375 & 0.060981 & $1.2\times10^{-2}$\tabularnewline
  \hline 
  Ex12, $A_{2}$, $\alpha=40$ & 0.04056702 & 0.03601258 & 0.060981 & $9.5\times10^{-3}$\tabularnewline
  \hline 
  Ex12, $A_{3}$, $\alpha=40$ & 0.06126412 & 0.06126407 & 0.060981 & $6.3\times10^{-7}$\tabularnewline
  \hline 
  Ex12, $A_{4}$, $\alpha=40$ & 0.06126593 & 0.06126569 & 0.060981 & $6.3\times10^{-7}$\tabularnewline
  \hline
\end{tabular}
\caption{Results on digital options for $d=4$\newline with
$r=0.05,T=1,S_0^{(1)}=S_0^{(2)}=S_0^{(3)}=S_0^{(4)}=50,U_1=U_2=U_3=U_4=60$
\newline
Ex11: $\rho=0.1,K=45$\newline
Ex12: $r=0.05,\rho=0.9,K=55$}\label{tab_4DDigital}
\end{table}
The number of function evaluations are $4.4\times10^{8}$ for $\alpha=30$ and
$5.8\times10^{8}$ for $\alpha=40$ meanwhile there are $2.2\times10^{7}$ and
$4.3\times10^{6}$ Monte Carlo runs respectively for Ex11 and Ex12.
Obviously, the complexity of the GRS increases with the dimension but it is
still more efficient on Ex11 than the Monte Carlo method.  On Ex12, the
volume of the subdomain where the function is positive is getting too small
compared to the size $A_{2}^4$ of the initial hypperrectangle for the GRS to
converge. Nevertheless, we have been able to recover a very good accuracy on
Ex12 by drastically reducing the size of the initial hypercube using $A_{3}$
and $A_{4}$. Indeed, the corresponding price values have almost 6 common
digits and the error criterion is $2\times10^{-7}.$ This shows that is would
be worth using some basic preliminary numerical methods to circumvent the
region of interest in a relatively small box.

\section{Delta computation}

\subsection{Tchebychef Interpolation and differentiation}
\label{sec-interpol}

To compute the partial derivatives involved in the computation of the Delta,
we use a standard approach based on Tchebychef interpolation polynomials. We
describe in detail the two dimensional case for a general function $f(x,y)$.
First, we compute the Tchebychef interpolation polynomial
$T_{m,h,x_0,y_0}(x)$ of $f(x,y_0)$ using $m$ interpolation points in the
interval $]x-h,x+h[$. The approximation of the first component of the
gradient of $f$ at point $(x_0,y_0)$ is given by $\displaystyle\frac{\partial
}{\partial x}T_{m,h,x_0,y_0}(x)$. Thus, it requires to compute $m$ prices,
which is done using the GRS. Using different values of the parameters $h$ and
$m$, we are able to obtain reliable Delta approximations.

\subsection{Monte Carlo approach} To measure the efficiency of the approach
described in Paragraph~\ref{sec-interpol}, we need an alternative method to
compute the Delta, which is sufficiently accurate to serve as a benchmark. We
relied on a finite difference approach coupled with a Monte Carlo technique.
We closely followed the recommendations of~\cite{key-GLYNN} to tune the
number of Monte Carlo simulation and the finite difference step accordingly
to ensure the best possible accuracy. For a Monte Carlo method with $n$
samples, we used a finite difference step equal to $h_n = n^{-1/6}$ so that
under mild assumptions the finite difference estimator converges almost
surely to the exact value and satisfies a central limit theorem with the rate
$n^{-1/3}$.

\subsection{Numerical Results}
For symmetry reasons, we compute only the first component of the Delta. 
We compare in table \ref{tab_3DDelta} the Monte Carlo approach and the method based on Tchebychef interpolation on 4 examples. These latter are taken in dimension 3 and 4 for basket and minimum options.
For the interpolation method, we use two sets of parameters chosen to have two different accurate estimations of the Delta. The maximum number of function evaluations required for the examples in dimension 4 is $5\times4.4\times10^{7}=2.2\times10^{8}$ which is comparable to the number of samples $n$ used in the Monte Carlo method.
\begin{table}[H]
\centering\begin{tabular}{|c|c|c|c|c|}
  \hline 
  $\ $ & Monte Carlo & $m=3,h=0.05$ & $m=5,h=0.1$ & $n$ \tabularnewline
  \hline 
  Ex13 & 0.300088 & 0.3002853 & 0.3002864  & $1.4\times10^{8}$\tabularnewline
  \hline 
  Ex14 & -0.240865 & -0.2382143 & -0.2382098  & $1.6\times10^{8}$\tabularnewline
  \hline 
  Ex15 & 0.230224 & 0.2303219 & 0.2303214 & $1.1\times10^{8}$\tabularnewline
  \hline 
  Ex16 & -0.186628 & -0.1837101 & -0.1836998 & $1.6\times10^{8}$\tabularnewline
   \hline
\end{tabular}
\caption{Delta Computations with $T=1, r=0.05, \sigma=0.2, S_0=50$
\newline
Ex13: Basket options, $d=3,\rho=0.1,K=45$\newline
Ex14: Put on minimum options, $d=3,\rho=0.5,K=55$\newline
Ex15: Basket options, $d=4,\rho=0.1,K=45$\newline
Ex16: Put on minimum options, $d=4,\rho=0.5,K=55$}\label{tab_3DDelta}
\end{table}
The value of $n$ is chosen to ensure an accuracy of 3 digits on the Delta.
All estimations of the Delta using the interpolation method have about 3 common digits with the Monte Carlo estimation. Since for the two sets of parameters we have obtained at least 5 common digits for the two different estimations, we may conclude that the interpolation method is more efficient than the Monte Carlo approach.

\section{Dimension reduction and control variates for high dimensional
problems} 
\subsection{Description of the method} 

The GRS developed and tested in the previous sections suffers from the curse
of dimensionality but shows a very impressive accuracy in low or medium
dimensions up to 4 or 5. In this section, we propose a way to use this method
to devise a control variate for high dimensional models. We are still
interested in computing expectations of the form $\E(\psi(S_{T}))$ where
$S_{T}$ is defined as in the previous section.  The basic idea of the method
is to perform a principal component analysis of $S_{T}$ in order to reduce
the dimension of the problem by keeping only the leading components and
setting the others to zero. The expectation in the reduced model can be computed
quickly and accurately using the GRS and can serve as a control variate for
the original problem.

\subsubsection{Principal Component Analysis} We rewrite the model to embed
the correlation matrix $\Gamma$ into the volatility structure turned into a
matrix $\Sigma$ defined by 
\[
\Sigma = \diag(\sigma_{1},...,\sigma_{d}) \; \Gamma \; 
\diag(\sigma_{1},...,\sigma_{d}).
\]
Note that the vector $(\sigma_{1}W_{T},...,\sigma_{2}W_{T})$ is a Gaussian
vector with covariance matrix $\sqrt{T} \Sigma$. It is straightforward to
check that the matrix $\Sigma$ inherits its symmetric positive definite
feature from the one of $\Gamma$ and hence admits an orthonormal basis of
eigenvectors with positive eigenvalues.  Let $D$ be the diagonal matrix built
up with these eigenvalues sorted in decreasing order and let $P$ be the
corresponding matrix of eigenvectors sorted accordingly.  Note that
reordering the columns of $P$ does not change its orthonormal property and
that we have $P^{-1}=P^{t}$ and as a consequence $\Sigma=P^{t}DP$. Now, if we
let $D^{\frac{1}{2}}=\diag(\sqrt{D_{11}},...,\sqrt{D_{dd}})$ and define the
symmetric matrix $H=P^{t}D^{\frac{1}{2}}P$, we obtain $\Sigma=H^{t}H$. This
leads to the identity in distribution
\[
(\sigma_{1}W_{T},...,\sigma_{2}W_{T})=\sqrt{T}HG
\]
where $G$ is a standard normal vector in $\mathbb{R}^{d}.$ Then, we can
write the identity in law
\[
S_{T}^{i}=S_{0}^{i}\exp\left((r-\frac{\sigma_{i}^{2}}{2})T+\sqrt{T}H_{i}G\right)
\]
where $H_{i}$ denotes the $i-{th}$ row of $H.$ This new expression for $S$
can help us to devise a reduced model $\hat{S}$ with effective dimension
$l\leq d$ such that $\hat S$ can be written
\[
\hat{S_{T}^{i}}=S_{0}^{i}\exp\left(r-\frac{\sigma_{i}^{2}}{2})T+\sqrt{T}H_{i}\hat{G}\right)
\]
with $\hat{G}=(G_{1},...,G_{l},0,...,0).$

\subsubsection{Control variates} If $l$ is sufficiently small, namely less
than 3 or 4, we can compute an approximation $\hat{I}$ of
$\E(\psi(\hat{S_{T}}))$ using the GRS both very accurately and with a small
computational cost.  If the dimension reduction works well, we can also
expect the random variable $\psi(\hat{S_{T}})$ to be close to $\psi(S_{T})$
and as a consequence to appear as a natural choice for a control variate.
Our new Monte Carlo estimator for $\E(\psi(S_{T}))$ will be 
\[
\frac{1}{N_{MC}}\sum_{j=1}^{N_{MC}}(\psi(S_{T}^{(j)})-\psi(\hat{S_{T}^{(j)}}))+\hat{I}
\]
where the $N_{MC}$ random samples $S_{T}^{(j)}$ and $\hat{S_{T}^{(j)}}$ are
built using the same independent drawings $G^{(j)}$ from the standard normal
distribution over $\mathbb{R}^{d}.$ Because the effective dimension $l$ of
$\hat{G}$ is a lot smaller than $d$ and the GRS is so efficient in low
dimensions, one can implement this control variate approach with almost no
extra computational cost. The idea of using such a control variate  was
already proposed for example in \cite{key-HENDERSON} but it was limited to a
sum of control variates depending on only one variable. Indeed this
computation was relying on closed forms for the expectations of the control
variates which are only available in dimension one. Hence, our new approach
is far more general.

\subsection{Numerical results}

\subsubsection{Toy correlation structures}

All the examples deal with basket options in dimension 5 with three different
values of the correlation coefficient $\rho$ from very correlated assets
($\rho=0.9$) to negatively correlated assets ($\rho=-0.1$). The different
estimations of the price and its confidence interval, named $I_{E_0}$,
$I_{E_1}$, $I_{E_2}$, $I_{E_3}$, using the control variate method are given
in table \ref{tab_5DACP2}. We compute the estimation of the price $I_5$ using
the GRS in order to obtain a reference value. 
\begin{table}[H]
\centering\begin{tabular}{|l|c|c|c|c|c|}
\hline 
  $\ $ & $I_{E_0}$  & $I_{E_1}$ & $I_{E_2}$ & $I_{E_3}$& $I_5$ \tabularnewline
\hline
Ex17&$8.61236\pm0.020$ &$8.61333\pm0.0013$ &$8.61357\pm0.0010$  &$8.61407\pm0.0003$ &$8.61404$ \tabularnewline
\hline 
Ex18&$7.51683\pm0.0130$ &$7.52217\pm0.0072$ &$7.52395\pm0.0042$ &$7.52621\pm0.0023$&$7.52490$ \tabularnewline
\hline 
Ex19&$7.29012\pm0.0103$ &$7.28436\pm0.0074$ &$7.27969\pm0.0042$ &$7.27931\pm0.0038$ &$7.27548$ \tabularnewline
\hline
\end{tabular}
\caption{Basket options with 0.95\% confidence interval\newline$d=5, T=1, r=0.05,K=45, S_0=50\newline\sigma=(0.156,0.442,0.325,0.134,0.114)$  and $100,000$ samples
\newline
Ex17: $\rho=0.9$, Ex18: $\rho=0.1$, Ex19: $\rho=-0.1$}
\label{tab_5DACP2}
\end{table}
We observe that the control variate technique gives very good results.
Especially with 3 PCA components, the results are very close to $I_5$.
Moreover, the GRS becoming very costly in dimension 5,  the control variate approach 
with 3 components seems now more efficient.

As it might be expected, we also noticed that when the assets are very
correlated the variance reduction is more efficient.  With $\rho=0.9$, the
first component of the PCA contains almost all the information of the model
which ensures a huge variance reduction. When $\rho=-0.1$, even with 3 PCA
components, the variance reduction is so small that it hardly compensates the
additional cost of the control. Indeed this latter is twice bigger than the
initial Monte Carlo cost.

\subsubsection{A more general correlation structure}

In this paragraph, we are concerned with testing our approach on more
realistic covariance structures. We consider a block covariance matrix
$G_{Ex20}$, in which the assets can be split into two subsets which are
negatively correlated, whereas inside each set all the assets are positively
correlated with the same correlation. Such covariance structures are quite
common in practice. We can find in Table~\ref{tab_10DACP} the results
obtained with a standard Monte Carlo and our control variate approach method
for different numbers of PCA components. At first sight, one could think that
such a covariance structure could cause trouble to our approach and yet it
manages to divide the width of the confidence interval by a factor of $10$,
which is really impressive if one remembers that to achieve such a confidence
interval with a crude Monte Carlo it would require $100$ times more samples.
We can see on this example that to effectively reduce the variance, we need
at least two PCA components and no closed form formulae are available for
basket options in dimension $2$ and $3$. This example highlights the key role
played by the GRS approach.

\[
\Gamma_{Ex20} = \left(\begin{array}{cccccccccc}
  1	& 0.8	& 0.8	& 0.8	& 0.8	& -0.5	& -0.5	& -0.5	& -0.5	& -0.5	\\
  0.8	& 1	& 0.8	& 0.8	& 0.8	& -0.5	& -0.5	& -0.5	& -0.5	& -0.5	\\
  0.8	& 0.8	& 1	& 0.8	& 0.8	& -0.5	& -0.5	& -0.5	& -0.5	& -0.5	\\
  0.8	& 0.8	& 0.8	& 1	& 0.8	& -0.5	& -0.5	& -0.5	& -0.5	& -0.5	\\
  0.8	& 0.8	& 0.8	& 0.8	& 1	& -0.5	& -0.5	& -0.5	& -0.5	& -0.5	\\
  -0.5	& -0.5	& -0.5	& -0.5	& -0.5	& 1	& 0.4	& 0.4	& 0.4	& 0.4	\\
  -0.5	& -0.5	& -0.5	& -0.5	& -0.5	& 0.4	& 1	& 0.4	& 0.4	& 0.4	\\
  -0.5	& -0.5	& -0.5	& -0.5	& -0.5	& 0.4	& 0.4	& 1	& 0.4	& 0.4	\\
  -0.5	& -0.5	& -0.5	& -0.5	& -0.5	& 0.4	& 0.4	& 0.4	& 1	& 0.4	\\
  -0.5	& -0.5	& -0.5	& -0.5	& -0.5	& 0.4	& 0.4	& 0.4	& 0.4	& 1	
\end{array}\right)
\]

\begin{table}[H]
  \centering\begin{tabular}{|l|c|c|c|c|}
    \hline 
    & $E_0$ & $E_1$ & $E_2$ & $E_3$ \tabularnewline
    \hline 
    Ex20 & $3.1912 \pm 0.011$ & $3.1899 \pm 0.009$ & $3.1908 \pm 0.002$ & $3.1906 \pm 0.001$\tabularnewline
    \hline 
  \end{tabular}
  \caption{Prices with 0.95\% confidence interval for homogeneous call Basket
  options on example Ex20 with
  \newline $T=2, r=0.02, K=105,\sigma=0.2, S_0=100$ and $100,000$ samples}
  \label{tab_10DACP}
\end{table}

\section{Conclusion}

In this paper, we have developed a new numerical integration algorithm based
on a Geometrical Random Splitting.  This new algorithm has been successfully
applied to the pricing of Vanilla options in  dimensions up to five.  In
particular, this new algorithm efficiently handles the difficult problems of
pricing and hedging digital options.  In most situations, the accuracies we
have obtained were out of reach for a crude Monte Carlo approach. For higher
dimensions, we have shown that the GRS can be used as a control variate when
associated to a principal component analysis.  The resulting variance
reduction was quite impressive especially for very correlated assets.  One
important issue would be the development of other dimension reduction
techniques that could be coupled with the GRS. The GRS itself can obviously
have many other applications for instance in computer experiment or in the
adaptive approximation of partial differential equations.  The GRS algorithm
proposed in this paper is very satisfactory on most examples studied. 



\begin{thebibliography}{References}
  \bibitem{key-ATANASSOV}E. I. ATANASSOV, I. T. DIMOV, A new optimal Monte
    Carlo method for calculating integral of smooth functions, Monte Carlo
    Methods and Appl., Vol. 5, No. 2, pp. 149-167, 1999.

  \bibitem{key-DEELSTRA}G. DEELSTRA, J. LIINEV, M. VANMAELE, Pricing of arithmetic
    basket options by conditioning, Insurance: Mathematics and Economics,
    vol 34, pp. 55-77, 2004.

  \bibitem{key-CDL}C. DE LUIGI, S. MAIRE, Adaptive integration and approximation
    over hyper-rectangular regions with applications to basket option
    pricing, Monte Carlo Methods and Applications, 16, (3-4), pp.265-282,
    2010. 

  \bibitem{key-HENDERSON} S.~M.~T. EHRLICHMAN and S.~G. HENDERSON.
    Adaptive control variates for pricing multi-dimensional American
    options, Journal of Computational Finance, 11(1), pp.65--91, 2007.

  \bibitem{key-ETORE}P. ETORE, G. FORT, B. JOURDAIN, E. MOULINES, On
    adaptive stratification, To appear in Annals of operations research.

  \bibitem{key-GLYNN} P.~W. GLYNN, Optimization of stochastic systems via
    simulation.  Proceedings of the 21st conference on Winter simulation, WSC
    '89, pages 90--105, New York, NY, USA, 1989. ACM.

  \bibitem{key-HELLUY}P. HELLUY, S. MAIRE, P. RAVEL, Int\' egration num\' erique
    d'ordre \' elev\' e de fonctions r\' eguli\` eres ou singuli\` eres sur un intervalle,
    CR. Acad. Sci. Paris, S\' er. 1, 327, pp. 843-848, 1998. 

  \bibitem{key-JOURDAIN}B. JOURDAIN, J. LELONG, Robust adaptive importance
    sampling for normal random vectors, Annals of applied probability,19
    (5), pp. 1687-1718, 2009.

  \bibitem{key-KROMMER}
    A. R. KROMMER, C. W. UEBERHUBER. Computational integration.
    SIAM, 1998.

  \bibitem{key-LAMBERTON}
    D. LAMBERTON, B. LAPEYRE, Introduction to stochastic calculus applied to finance,
    Chapman \& Hall, London, 1996.

  \bibitem{key-LAPEYRE}
    B. LAPEYRE and J. LELONG,
    A framework for adaptive Monte-Carlo procedures,
    Monte Carlo Methods Appl., vol. 17 (1), pp.77-98, 2011.

  \bibitem{key-LEPAGE}G. P. LEPAGE, A New Algorithm for Adaptative Multidimensional
    Integration, Journal of Computational Physics, Vol 27, pp. 192-203,
    1978.

  \bibitem{key-MAIRE1}S. MAIRE, Reducing variance using iterated control variates,
    The Journal of Statistical Computation and Simulation, Vol. 73(1),
    pp. 1-29, 2003. 

  \bibitem{key-MAIRE2}S. MAIRE, An iterative computation of approximations
    on Korobov-like spaces. Journal of Computational and Applied Mathematics,
    157, pp. 261-281, 2003.

  \bibitem{key-MAIRE3}S. MAIRE, Polynomial Approximations of multivariate
    smooth functions from quasi-random data, Statistics and Computing,
    14, pp. 333-336, 2004.

  \bibitem{key-MAIRE4}S. MAIRE, C. DE LUIGI, Quasi-Monte Carlo quadratures
    for multivariate smooth functions, Applied Numerical Mathematics,
    56, no.2, pp. 146-162, 2006.

  \bibitem{key-NIEDERREITER}H. NIEDERREITER, Quasi-Monte Carlo methods and pseudorandom
    numbers, Bull. Amer. Math. Soc. 84, pp. 957-1041, 1978.

  \bibitem{key-NOVAK}E. NOVAK, K. RITTER, High dimensional integration
    of smooth functions over cubes, Numerishe Mathematik, 75, pp.79-97,
    1996.

  \bibitem{key-PAGES}G. PAGES, A space vector quantization for numerical
    Integration, Journal of computational and applied mathematics, 89,
    pp. 1-38, 1997.

  \bibitem{key-PRESS}W. H. PRESS, G. R. FARRAR, Recursive Stratified Sampling
    for Multidimensional Monte Carlo Integration, Computer in Physics,
    vol. 4, pp. 190-195, 1990.

  \bibitem{key-SCHURER1}R. SCHURER, Adaptive quasi-Monte Carlo integration
    based on MISER and VEGAS. Monte Carlo and quasi-Monte Carlo methods
    2002, pp. 393-406, Springer, Berlin, 2004.

  \bibitem{key-SCHURER2}R. SCHURER, A comparison between (quasi)-Monte Carlo
    and cubature rule based methods for solving high-dimensional integration
    problems, Mathematics and computers in simulation, 62, pp. 509-517,
    2003.

  \bibitem{key-SLOAN}I. H. SLOAN, P. J. KACHOYAN, Lattice methods for multiple
    integration: Theory, error analysis and examples, SIAM J. Numer. Anal.
    24, pp. 116-128, 1987.\end{thebibliography}
    \end{document}